\documentclass[12pt,a4paper]{amsart}

\usepackage {amsfonts}
\usepackage[english]{babel}
\def\g{\gamma}
\def\G{\Gamma}
\def\d{\delta}
\def\a{\alpha}
\def\b{\beta}
\def\p{\varphi}
\def\e{\varepsilon}
\def\l{\lambda}
\def\L{\Lambda}
\def\s{\sigma}
\def\t{\theta}
\def\R{{\mathbb R}}
\def\C{{\mathbb C}}
\def\N{{\mathbb N}}
\def\Z{{\mathbb Z}}

\def\o{\omega}
\def\p{\varphi}
\def\bs{~\hfill\rule{7pt}{7pt}}
\def\la{\langle}
\def\ra{\rangle}

\DeclareMathOperator{\supp}{supp}
\DeclareMathOperator{\ord}{ord}

\newtheorem{Th}{Theorem}
\newtheorem{Cor}{Corollary}
\newtheorem{Lem}{Lemma}
\newtheorem{Pro}{Proposition}
\newtheorem*{Def}{Definition}

\begin{document}

\title{Almost periodic distributions and crystalline measures}

\author{S.Yu. Favorov}

\address{Serhii Favorov,
\newline\hphantom{iii}Faculty of Mathematics and Computer Science, Jagiellonian University,
\newline\hphantom{iii} Lojasiewicza 6, 30-348 Krakow, Poland
\newline\hphantom{iii} Faculty of Mathematics and Informatics, Karazin's Kharkiv National University
\newline\hphantom{iii} Svobody sq., 4, Kharkiv, Ukraine 61022}
\email{sfavorov@gmail.com}

\maketitle {\small
\begin{quote}
\noindent{\bf Abstract.}
Based on the properties of distributions and measures with discrete support, we investigate temperate almost periodic distributions  on the Euclidean space
and connection with their Fourier transforms. We also study relations between the Fourier transform of almost periodic distributions and their Fourier coefficients.
The main result of the article is the construction of a crystalline measure on the real line, which is neither
almost periodic distribution, nor a Fourier quasicrystal.

\medskip


\end{quote}
}

\bigskip
The concept of  Fourier quasicrystal was inspired by experimental discovery of nonperiodic
atomic structures with diffraction patterns consisting of spots, made in the middle of 80's.
A number of papers has appeared, in which the properties of Fourier quasicrystals are studied.  Conditions for support of  quasicrystals to be a finite union
of discrete lattices are found, and nontrivial examples of quasicrystals  are constructed (\cite{Q},\cite{F1},\cite{F3}-\cite{F5},\cite{K1}-\cite{L1},
\cite{LO1}-\cite{LO4},\cite{M1}-\cite{M4}).
These studies have been extended to the more general setting of temperate distributions with discrete support and spectrum (\cite{LR},\cite{F2},\cite{F6},\cite{F7},
\cite{P}). Note that the properties of almost periodic measures and distributions were more or less explicitly used in these investigates.

The structure of this article is as follows. In Section \ref{S1} we give the necessary definitions and notations. In Section \ref{S2} we prove some properties
of distributions with locally finite support and spectrum. These results are a slight enhancement of the results of \cite{F2} and are included here for the sake
 of completeness. In section \ref{S3} we introduce a notion of s-almost periodicity of temperate distributions, which is very close to the classical  almost periodicity.
 Here we present some theorems on connections between temperate distributions and their Fourier transform, some of them refer to almost periodic distributions,
 and others to s-almost periodic distributions. Based on them, we obtain the main result of the paper: the example of a crystalline measure that is not a Fourier quasicrystal.
 In Section \ref{S4} we give proofs of some theorems from Section \ref{S3}. Finally, in Section \ref{S5} we show that Meyer's Theorem  on the connection between
 the Fourier transform of almost periodic functions and measures and their Fourier coefficients \cite{M4} can be generalized to temperate almost periodic distributions
 and their Fourier coefficients in the sense of L.Ronkin \cite{R}.

          \section{Definitions and notations}\label{S1}
          \medskip

Denote by $S(\R^d)$ the Schwartz space of test functions $\p\in C^\infty(\R^d)$ with the finite norms
 $$
  N_{n,m}(\p)=\sup_{\R^d}\{\max\{1,|x|^n\}\max_{\|k\|\le m} |D^k\p(x)|\},\quad n,m=0,1,2,\dots,
 $$
where
$$
k=(k_1,\dots,k_d)\in(\N\cup\{0\})^d,\ \|k\|=k_1+\dots+k_d,\  D^k=\partial^{k_1}_{x_1}\dots\partial^{k_d}_{x_d}.
 $$
 These norms generate the topology on $S(\R^d)$.  Elements of the space $S^*(\R^d)$ of continuous linear functionals on $S(\R^d)$ are called temperate distributions.
 For each temperate distribution $f$ there are $C<\infty$ and $n,\,m\in\N\cup\{0\}$ such that for all $\p\in S(\R^d)$
\begin{equation}\label{d}
                           |f(\p)|\le CN_{n,m}(\p).
\end{equation}
Moreover, this estimate is sufficient for the distribution $f$ to belong to $S^*(\R^d)$
(see \cite{V}, Ch.3).

The Fourier transform of a temperate distribution $f$ is defined by the equality
$$
\hat f(\p)=f(\hat\p)\quad\mbox{for all}\quad\p\in S(\R^d),
$$
where
$$
   \hat\p(y)=\int_{\R^d}\p(x)\exp\{-2\pi i\la x,y\ra\}dx
 $$
is the Fourier transform of the function $\p$. By $\check\p$ we  denote the inverse Fourier transform of $\p$. The Fourier transform
is the bijection of $S(\R^d)$ on itself and the bijection of $S^*(\R^d)$ on itself.
The support of $\hat f$ is called {\it spectrum} of $f$.

  We will say that a set $A\subset\R^d$ is {\it locally finite} if the intersection of $A$ with any ball is finite,  $A$ is {\it relatively dense} if there is $R<\infty$
  such that $A$ intersects with each ball of radius $R$, and $A$ is {\it uniformly discrete}, if $A$ is locally finite and has a strictly positive separating constant
 $$
 \eta(A):=\inf\{|x-x'|:\,x,\,x'\in A,\,x\neq x'\}.
  $$
  Also, we will say that $A$ is {\it polynomially discrete}, or shortly {\it p-discrete}, if there are positive numbers $c, h$ such that
  \begin{equation}\label{c}
    |x-x'|\ge c\min\{1,\,|x|^{-h},|x'|^{-h}\}\qquad \forall x, x'\in A,\quad x\neq x'.
  \end{equation}
    A  set $A$ is of {\it bounded density} if it is locally finite and
  $$
        \sup_{x\in\R^d}\# A\cap B(x,1)<\infty.
  $$
  As usual, $\# E$ is a number of elements of the finite set $E$, and $B(x,r)$ is the ball with center at the point $x$ and radius $r$.

  An element $f\in S^*(\R^d)$ is called {\it a crystalline measure} if $f$ and $\hat f$ are  complex-valued measures on $\R^d$ with locally finite supports.

  Denote by $|\mu|(A)$ the variation of the complex-valued measure $\mu$ on $A$. If both measures $|\mu|$ and $|\hat\mu|$ have  locally finite supports and belong to $S^*(\R^d)$,
we say that $\mu$ is a {\it Fourier quasicrystal}. A measure $\mu=\sum_{\l\in\L}a_\l\d_\l$ with $a_\l\in\C$ and countable $\L$ is called {\it purely point},
here $\d_y$ is  the unit mass at the point $y\in\R^d$. If this is the case, we will replace $a_\l$ with $\mu(\l)$ and write $\supp\mu=\{\l:\,\mu(\l)\neq0\}$.

     \section{Temperate distributions with locally finite support}\label{S2}

\medskip

By \cite{Ru}, every distribution $f$ with locally finite support $\L$ has the form
$$
f=\sum_{\l\in\L} P_\l(D)\d_\l, \quad P_\l(x)=\sum_{\|k\|\le K_\l}p_k(\l)x^k, \quad x\in\R^d,\ p_k(\l)\in\C,\ K_\l<\infty.
$$
 Here, as usual,  $x^k=x_1^{k_1}\cdots x_d^{k_d}$. Note that $\ord f=\sup_\l\deg P_\l(x)\le\infty$.
.
\begin{Pro}\label{P1}
 Suppose $f\in S^*(\R^d)$ has a  locally finite support $\L$. Then

i) $\ord f<\infty$, hence,
\begin{equation}\label{r1}
f=\sum_{\l\in\L}\sum_{\|k\|\le K}p_k(\l)D^k\d_\l,\quad k\in(\N\cup\{0\})^d,\quad K=\ord f;
\end{equation}
in particular, if $f$ has a locally finite spectrum $\G$, then $\ord\hat f<\infty$ and
\begin{equation}\label{r2}
\hat f=\sum_{\g\in\G}\sum_{\|j\|\le J}q_j(\g)D^j\d_\g,\quad j\in(\N\cup\{0\})^d,\quad J=\ord\hat f.
\end{equation}

ii) If $\L$ is $p$-discrete,  then there exist $C,\,T<\infty$ such that for all $k$
 \begin{equation}\label{p}
  |p_k(\l)|\le C\max\{1,|\l|^T\} \quad\mbox{for all}\ \l\in\L.
 \end{equation}
Moreover, there exists $T_1<\infty$ such that
\begin{equation}\label{k}
\sum_{\l\in\L,|\l|<R}\sum_{\|k\|\le K}|p_k(\l)|=O(R^{T_1})\quad\mbox{as}\quad R\to\infty.
\end{equation}
   \end{Pro}

{\bf Proof of Proposition \ref{P1}}. i) Let $\l\in\L$ and  $\e\in(0,1)$ be such that
$$
\inf\{|\l-\l'|:\,\l'\in\L,\,\l'\neq\l\}>\e.
 $$
 Let $\p$ be a non-negative function on $\R$ such that
 \begin{equation}\label{p0}
 \p(|x|)\in C^\infty(\R^d),\quad\p(|x|)=0\mbox{ for }|x|>1/2,\quad \p(|x|)=1\mbox{ for }|x|\le1/3.
  \end{equation}
  Then set
  $$
  \p_{\l,k,\e}(x)=\frac{(x-\l)^k}{k!}\p\left(\frac{|x-\l|}{\e}\right)\in S(\R^d),
  $$
   where, as usual, $k!=k_1!\cdots k_d!$. It is easily shown that
   $$
   f(\p_{\l,k,\e})=(-1)^{\|k\|}p_k(\l).
   $$
Let $f$ satisfy (\ref{d}) with some $m,\,n$.  We get
$$
 |f(\p_{\l,k,\e})|\le C\sup_{|x-\l|<\e}\max\{1,|x|^n\}\sum_{\|\a+\b\|\le m} c(\a,\b)\left|D^\a\p\left(\frac{|x-\l|}{\e}\right)D^\b\left(\frac{(x-\l)^k}{k!}\right)\right|,
$$
 where $\a,\b\in(\N\cup\{0\})^d$ and $c(\a,\b)<\infty$. Note that
$$
\left|D^\a\p\left(\frac{|x-\l|}{\e}\right)\right|\le \e^{-\|\a\|}c(\a)\quad\mbox{for}\quad |\l-x|<\e/3,
$$
 and this derivative  vanishes for $|\l-x|\ge \e/2$. Also,
$$
D^\b(x-\l)^k=\begin{cases} 0 &\text{ if }k_j<\b_j\text{ for at least one }j,\\c(k,\b)(x-\l)^{k-\b} &\text{ if }k_j\ge\b_j\,\ \forall j.
  \end{cases}.
$$
Since
$$
\max\{1,|x|^n\}\le 2^n\max\{1,|\l|^n\}
$$
for $x\in\supp\p_{\l,k,\e}$, we get
$$
 |p_k(\l)|\le\sum_{\|\a+\b\|\le m,\,\b_j\le k_j\,\forall j} c(k,\a,\b)\max\{1,|\l|^n\}\e^{\|k\|-\|\a+\b\|}.
$$
For $\|k\|>m$  we take $\e\to0$ and obtain $p_k(\l)=0$.

Since $\hat f\in S^*(\R^d)$, we obtain \eqref{r2}.
\smallskip

ii) Let $\L$ be $p$-discrete and (\ref{p})  be not satisfy.  Then there is $k, \|k\|\le K,$ and a sequence $\l_s\to\infty$ such that $|\l_{s+1}|>1+|\l_s|$ for all $s$ and
\begin{equation}\label{p4}
  \log|p_k(\l_s)|/\log|\l_s|\to\infty, \quad s\to\infty.
 \end{equation}
 Put $\b_s=c|\l_s|^{-h}$ with $c$ from \eqref{c} and
 $$
 \psi_{s,k}(x)=\frac{(x-\l_s)^k}{k!}\p\left(\frac{|x-\l_s|}{\b_s}\right),\qquad
 \Psi_k(x)=\sum_{s=1}^\infty \frac{\psi_{s,k}(x)}{p_k(\l_s)}.
  $$
 We may suppose that $\l_1$ is large enough such that $\supp\psi_{s,k}\cap\supp\psi_{s',k}=\emptyset,\ s\neq s'$.  Then by (\ref{p4}),
$$
1/p_k(\l_s)=o(1/|\l_s|^T),\quad  |\l_s|\to\infty, \quad\mbox{ for every }T<\infty.
$$
 Since
 $$
 D^j(\psi_{s,k}(x))=O(|\l_s|^{h\|j\|}),\ j\in(\N\cup\{0\})^d,
 $$
  we see that
 $$
 D^j(\Psi_k(x))=o(1/|x|^{T-h\|j\|}),\quad x\to\infty,
 $$
 and $\Psi_k\in S(\R^d)$.

 Since $\L$ is $p$-discrete, we get $\l\not\in B(\l_s,c|\l_s|^{-h})$ for all $\l\in\L\setminus\{\l_s\}$. Therefore, $f(\Psi_k)$ is equal to
$$
   \sum_{\l\in\L}\sum_{\|l\|\le K}\sum_s(-1)^{\|l\|}p_l(\l)p_k(\l_s)^{-1}D^l(\psi_{s,k})(\l)=
   \sum_s\sum_{\|l\|\le K} (-1)^{\|l\|}p_l(\l_s)p_k(\l_s)^{-1}D^l(\psi_{s,k})(\l_s).
$$
 Since $D^l(\psi_{s,k})(\l_s)=0$ for $l\neq k$ and $D^k(\psi_{s, k})(\l_s)=1$, we obtain the contradiction.
\medskip

Estimate \eqref{k} follows immediately from \eqref{p} and  the following simple lemma:
\begin{Lem}[cf. \cite{F7}, a part of the proof of Theorem 6]\label{L1}
If $S$ is $p$-discrete set, then  $\# S\cap B(0,R)=O(R^{T'})$ as $R\to\infty$ with $T'<\infty$.
\end{Lem} \bs

{\bf Remark}.  Proposition \ref{P1} was earlier proved by V.Palamodov \cite{P} for temperate distributions with uniformly discrete support.

\begin{Pro}\label{P2}
Let $\mu\in S^*(\R^d)$ be a measure. Then $|\mu|$ belongs to $S^*(\R^d)$ if and only if there is $T<\infty$ such that $|\mu|(B(0,R))=O(R^T)$ as $R\to\infty$.
 \end{Pro}

{\bf Proof}. Any non-negative measure $\nu$ on $\R^d$ satisfying the  condition $\nu(B(0,R))=O(R^T)$ as $R\to\infty$ belongs to $S^*(\R^d)$ (cf.\cite{Ru}).
The converse statement see \cite{F5}, Lemma 1. \bs
\medskip

It follows from Propositions \ref{P1} and \ref{P2}
\begin{Th}[cf. also \cite{F7}]\label{T1}
Let a measure $\mu$ has p-discrete support and belongs to $S^*(\R^d)$. Then $|\mu|\in S^*(\R^d)$ too. In particular, every crystalline measure with p-discrete support
and p-discrete spectrum is the Fourier quasicrystal.
\end{Th}

 M.Kolountzakis, J.Lagarias proved in \cite{KL} that the Fourier transform  of every measure $\mu$ on the line $\R$ with locally finite support of bounded density,
 bounded masses $\mu(x)$, and locally finite spectrum is also a measure $\hat\mu=\sum_{\g\in\G}q_\g\d_\g$ with uniformly bounded $q_\g$.
 The following proposition generalizes this result for distributions from $S^*(\R^d)$.

 \begin{Pro}\label{P3}
Suppose $f\in S^*(\R^d)$  has form (\ref{r1}) with some $K$ and countable $\L$, and $\hat f$ has form \eqref{r2} with the locally finite support $\G$. If
 $$
\rho_f(r):=\sum_{|\l|<r}\sum_{\|k\|\le K}|p_k(\l)|=O(r^{d+H}),\quad r\to\infty,\quad H\ge0,
  $$
  then $\ord\hat f\le H$; if $\rho_f(r)=o(r^{d+H})$ as $r\to\infty$, then $\ord\hat f<H$.

Furthermore, in the case of integer $H$ and $\|j\|=H$ we have $|q_j(\g)|\le C'\max\{1,|\g|^K\}$;
for the case of uniformly discrete $\G$ this estimate with the same $K$ takes place for all $j$.
\end{Pro}
\begin{Cor}\label{C1}
 If $f\in S^*(\R^d)$  has form (\ref{r1}) with countable $\L$, locally finite spectrum $\G$,  and $\rho_f(r)=O(r^d)$ as $r\to\infty$,
 then $\hat f$ is a measure, and
 $$
 \hat f=\sum_{\g\in\G}q(\g)\d_\g,\ |q(\g)|\le C'\max\{1,|\g|^K\}.
 $$
\end{Cor}

 {\bf Proof of Proposition \ref{P3}}. Let $\g\in\G$ and pick $\e\in(0,1)$ such that
$$
\inf\{|\g-\g'|:\,\g'\in\G,\,\g'\neq\g\}>\e.
 $$
 Let $\p$ be the same as in the proof of Proposition \ref{P1}. Put
 $$
 \p_{\g,l,\e}(y)=\frac{(y-\g)^l}{l!}\p(|y-\g|/\e)\in S(\R^d).
 $$
  We have
 $$
 (-1)^{\|l\|}q_l(\g)=\sum_{\|j\|\le J}q_j(\g)D^j\d_\g(\p_{\g,l,\e}(y))=(\hat f,\p_{\g,l,\e})=(f,\hat\p_{\g,l,\e}).
 $$
Note that
 $$
 \hat\p_{\g,l,\e}(x)=e^{-2\pi i\la x,\g\ra}(l!)^{-1}(-2\pi i)^{-\|l\|}D^l(\widehat{\p(\cdot/\e)})=c(l)e^{-2\pi i\la x,\g\ra}\e^{d+\|l\|}(D^l\hat\p)(\e x).
 $$
 Therefore,
 $$
 D^k(\hat\p_{\g,l,\e})(x)=\e^{d+\|l\|}\sum_{\a+\b=k}c(\a,\b)D^\a\left[e^{-2\pi i\la x,\g\ra}\right]D^\b[(D^l\hat\p)(\e x)]
 $$
 $$
 =\sum_{\a+\b=k}c(\a,\b)(-2\pi i)^{\|\a\|}\g^\a e^{-2\pi i\la x,\g\ra}\e^{d+\|l\|+\|\b\|}(D^{\b+l}\hat\p)(\e x).
 $$
 Since $\hat\p(\e x)\in S(\R^d)$, we get for every $x\in\R^d$ and $n\in\N\cup\{0\}$
 $$
   |D^{\b+l}(\hat\p)(\e x)|\le N_{n,\|\b+l\|}(\hat\p)(\max\{1,|\e x|^n\})^{-1}.
 $$
 Therefore for every $k$, $\|k\|\le K$,
  $$
 |D^k(\hat\p_{\g,l,\e})(x)|\le C(K,n)\e^{d+\|l\|}\max\{1,|\g|^K\}(\max\{1,|\e x|^n\})^{-1},
 $$
where $C(k,n)$ depends on $\p$. Now we may estimate $(f,\hat\p_{\g,l,\e})$ as
\begin{equation}\label{h}
  \left|\sum_k\sum_\l p_k(\l)D^k(\hat\p_{\g,l,\e})(\l)\right|\le
 C(K,n)\e^{d+\|l\|}\max\{1,|\g|^K\}\int_0^\infty\frac{\rho_f(dt)}{\max\{1;\,(\e t)^n\}}.
 \end{equation}
If $\rho_f(r)=O(r^{d+H})$ for $r\to\infty$, take $t_0$ such that $|\rho(t)|<C_0t^{d+H}$ for $t>t_0$. If $\rho_f(r)=o(r^{d+H})$, fix any $\eta>0$ and take
$t_0=t_0(\eta)$ such  that  $|\rho(t)|<\eta t^{d+H}$ for $t>t_0$. Then pick $n>d+H$ and $\e<1/t_0$. Integrating by parts and using the estimate for $\rho_f(t)$, we obtain
$$
 \int_0^\infty\max\{1,(\e t)^n\}^{-1}\rho_f(dt)= \rho_f(1/\e)+\int_{1/\e}^\infty(\e t)^{-n}\rho_f(dt)\le \frac{nC_0}{\e^n}\int_{1/\e}^\infty t^{d+H-n-1}dt.
$$
Therefore, the left-hand side of \eqref{h} not more than  $\e^{\|l\|-H} C_0C'\max\{1,|\g|^K\}$, and
  $$
 |q_l(\g)|\le C'C_0\max\{1,|\g|^K\}\e^{\|l\|-H}.
 $$

  If $\|l\|>H$, we take $\e\to0$ and get $q_l(\g)=0$, hence, $J=\ord\hat f\le H$.
\smallskip

 If $H$ is integer, we get $|q_l(\g)|\le C'C_0\max\{1,|\g|^K\}$ for $\|l\|=H$.
\smallskip

 If $\rho_f(r)=o(r^{d+H})$, we replace $C_0$ by $\eta$ and note that $\eta$ is arbitrary small for $\e$ small enough. Hence,
$q_l(\g)=0$ for $\|l\|=H$.

 Finally, if $\G$ is uniformly discrete,  we take $\e=\e_0<\eta(\G)/2$ for all $\g\in\G$ and obtain the bound
  $$
 |q_l(\g)|\le \e_0^{-H} C'C_0\max\{1,|\g|^K\}\quad \forall l, \|l\|\le J.
 $$

 \bs

\section{Almost periodic distributions and their properties}\label{S3}
\medskip

Recall that a continuous function $g$ on $\R^d$ is  almost periodic if for any  $\e>0$ the set of $\e$-almost periods of $g$
  $$
  \{\tau\in\R^d:\,\sup_{x\in\R^d}|g(x+\tau)-g(x)|<\e\}
  $$
is a relatively dense set in $\R^d$.

Almost periodic functions are uniformly bounded  on $\R^d$. The class of almost periodic functions is closed with respect to taking absolute values,
 and finite linear combinations;  the limit of a uniformly in $\R^d$ convergent sequence  of almost periodic functions is also almost periodic.

   A typical example of an almost periodic function is every absolutely convergent exponential
   sum $\sum c_n\exp\{2\pi i\la x,\o_n\ra\}$ with $\o_n\in\R^d,\,c_n\in\C$ (cf.,for example, \cite{C}, \cite{M4}).

A measure $\mu$ on $\R^d$ is called almost periodic if the function
 $$
 (\psi\star\mu)(t)=\int_{\R^d}\psi(t-x)d\mu(x)
 $$
  is almost periodic in $t\in\R^d$ for each continuous function $\psi$ on $\R^d$ with compact support. A distribution $f\in S^*(\R^d)$ is almost periodic
  if the function $(\psi\star f)(t)=f(\psi(t-\cdot))$ is almost periodic in $t\in\R^d$ for each $\psi\in C^\infty$ with compact support
   (see \cite{LA}, \cite{R}, \cite{M1}, \cite{M2}, \cite{M4}). Clearly, every almost periodic distribution has a relatively dense support.
    But there are measures that are almost periodic temperate distributions, but are not almost periodic as measures (see \cite{M1}).

 \begin{Def}
  A distribution $f\in S^*(\R^d)$ is s-almost periodic, if the function $(\psi\star f)(t)=f(\psi(t-\cdot))$ is almost periodic in $t\in\R^d$ for each $\psi\in S(\R^d)$.
\end{Def}
The following theorem plays a very important role in our investigations.
\begin{Th}\label{T2}
If $f$ is a temperate distribution and its Fourier transform $\hat f$ is a pure point measure such that $|\hat f|(B(0,r))=O(r^T)$ for $r\to\infty$  with some $T<\infty$,
then $f$ is s-almost periodic  distribution.
\end{Th}
The proof is very easy. Let $\hat f=\sum_{\g\in\G} b(\g)\d_\g$, $M(r)=|\hat f|(B(0,r))$. For any $\psi\in S(\R^d)$ we have
\begin{equation}\label{d1}
  f(\psi(t-\cdot))=(\hat f(y),\hat\psi(y)e^{2\pi i\la t,y\ra})=\sum_{\g\in\G} b(\g)\hat\psi(\g)e^{2\pi i\la t,\g\ra}.
\end{equation}
 Since $|\hat\psi(y)|\le N_{T-1,0}(\hat\psi)|y|^{-T-1}$ for $|y|>1$ and
$$
  \sum_{\g\in\G} |b(\g)||\hat\psi(\g)|\le C_0+C_1\int_1^\infty r^{-T-1}M(dr)<\infty,
$$
we see that the series in \eqref{d1} absolutely converges, and the function $(f\star\psi)(t)$ is almost periodic.  \bs
\smallskip

 From Proposition \ref{P2} it follows
\begin{Cor}\label{C2}
If $f\in S^*(\R^d)$, $\hat f$ is pure point measure, and $|\hat f|\in S^*(\R^d)$, then $f$ is s-almost periodic  distribution.
 In particular, every Fourier quasicrystal is s-almost periodic  distribution.
\end{Cor}
\medskip

Using Proposition \ref{P1}, we also get
\begin{Cor}\label{C3}
  Let $f\in S^*(\R^d)$  have $p$-discrete spectrum, and the Fourier transform $\hat f$ is a measure. Then $f$ is s-almost periodic  distribution.
 \end{Cor}

\begin{Cor}\label{C4}
 Let $f\in S^*(\R^d)$ of form \eqref{r1} have locally finite  spectrum $\G$ with polynomial growth of numbers $\#(\G\cap B(0,r))$. If
$$
   \sum_{\|k\|\le K}\sum_{|\l|<r}|p_k(\l)|=O(r^d) \quad\mbox{for}\quad r\to\infty,
$$
 then $f$ is s-almost periodic  distribution.
  \end{Cor}
  Indeed, by Corollary \ref{C1}, $\hat f$ is a measure with polynomially growth coefficients.
\medskip

Evidently, every s-almost periodic  distribution is an almost periodic distribution too. I do not know if there is an almost periodic distribution,
 which is not s-almost periodic. However, the following assertion is valid:
\begin{Th}\label{T3}
  Every almost periodic (in sense of distributions) non-negative measure $\mu\in S^*(\R^d)$ is s-almost periodic distribution. The same implication is valid
  if $\mu$ is a complex-valued measure on $\R^d$ such that
\begin{equation}\label{a}
  \sup_{x\in\R^d}|\mu|(B(x,1))<\infty.
 \end{equation}
\end{Th}
It is easy to check that every almost periodic in the sense of distributions measure $\mu$ under condition \eqref{a} is almost periodic
in sense of measures.

\medskip
Proofs of Theorem \ref{T3} and the following ones are given in the next Section \ref{S4}.

\begin{Th}\label{T4}
  If $f\in S^*(\R^d)$ is an almost periodic distribution with locally finite spectrum $\G$, then $\hat f$ is a measure.
  \end{Th}

Show that $p$-discreteness of support of a measure is closely connected with s-almost periodicity  of its Fourier transform:
\begin{Th}\label{T5}
In order for each measure $\mu\in S^*(\R^d)$ with support in a fixed locally finite set $A\subset\R^d$ to have s-almost periodic  Fourier transform $\hat\mu$,
it is necessary and sufficient that $A$ be $p$-discrete.

Moreover, if $\hat\mu\star\psi(t)$ is bounded for all $\psi\in S(\R^d)$ and $\mu\in S^*(\R^d)$ with $\supp\mu\subset A$, then $A$ is $p$-discrete too.
 \end{Th}

\begin{Th}\label{T6}
There is a crystalline measure $\mu$ on $\R$ such that for some $\psi\in S(\R)$ the function $(\mu\star\psi)(t)$ is
unbounded in $t\in\R$. In particular,  $\mu$ is neither s-almost periodic  distribution, nor the Fourier quasicrystal.
\end{Th}

{\bf Remark}. Y.Meyer formulated in \cite{M2} as a theorem that any crystalline measure is an almost periodic distribution. Then he wrote in \cite{M3} that
the proof of this theorem is incorrect and formulated the corresponding result as Conjecture 2.1.  Theorem \ref{T6} gives only a partial answer on this conjecture,
because the function $\psi\in S(\R)$, which constructs in the proof of Theorem \ref{T6}, does not have compact support. We don't know if there is a function
 $\tilde\psi\in C^\infty$ with compact support  such that $(\mu\star\tilde\psi)(t)$ is unbounded.

\section{Proofs of the theorems}\label{S4}
\medskip

{\bf Proof of Theorem \ref{T3}}. Let $\p$ be $C^\infty$ non-negative function with compact support such that $\p(x)\equiv1$ for $x\in B(0,1)$. Since $\p\star\mu(t)$ is an almost
periodic function, we see that it is uniformly bounded. If $\mu\ge0$, we get $\mu(B(x,1))<C$ for all $x\in\R^d$.

Set $\mu^t(x):=\mu(t-x)$ with $t\in\R^d$. For every complex-valued measure $\mu$ under condition
\eqref{a} we get $M(r):=|\mu^t|(B(0,r))<Cr^d$ for all $r>1$, where the constant $C$ is the same for all $t$.
Take $\psi\in S(\R^d)$. Then $|\psi(x)|\le C_1|x|^{-d-1}$ for $|x|>1$. For any $\e>0$ there is $R<\infty$ that does not depend on $t$ and such that
$$
  \left|\int_{|x|>R}\psi(x)\mu^t(dx)\right|\le C_1\int_R^\infty r^{-d-1}M(dr)\le C_1(d+1)\int_R^\infty M(r)r^{-d-2}dr<\e/3.
$$
Therefore for all $t\in\R^d$
\begin{equation}\label{b}
  \left|\int_{|t-x|>R}\psi(t-x)\mu(dx)\right|=\left|\int_{|x'|>R}\psi(x')\mu^t(dx')\right|<\e/3.
 \end{equation}
Let $\xi(x)$ be $C^\infty$-function on $\R^d$ such that
$$
0\le\xi\le1,\quad \xi(x)\equiv1\quad\text{for}\quad |x|<R,\quad \xi(x)\equiv0\quad\text{for}\quad |x|>R+1.
$$
The function $(\xi\psi)\star\mu(t)$ is almost periodic, hence there are a relatively dense set $E\subset\R^d$ such that for any $\tau\in E$ and all $t\in\R^d$
$$
|(\xi\psi)\star\mu(t+\tau)-(\xi\psi)\star\mu(t)|<\e/3.
$$
Applying \eqref{b} to $(1-\xi(t+\tau))\psi(t+\tau)$ and $(1-\xi(t))\psi(t)$, we obtain
$$
|\psi\star\mu(t+\tau)-\psi\star\mu(t)|\le|(\xi\psi)\star\mu(t+\tau)-(\xi\psi)\star\mu(t)|+|(1-\xi)\psi\star\mu(t+\tau)|+|(1-\xi)\psi\star\mu(t)|<\e.
$$
Hence, $E$ is the set of $\e$-almost periods for the function $\psi\star\mu$.  \bs

\medskip

{\bf Proof of Theorem \ref{T4}}.
  Let $f$ be an almost periodic temperate distribution with a locally finite spectrum $\G$. By Proposition \ref{P1}, $\hat f$ has  form (\ref{r2}).
  Suppose that $J\neq0$ and $q_{j'}(\g')\neq0$ for some $\g'\in\G$ and $j'=(j'_1,\dots,j'_d),\,\|j'\|=J$. Without loss of generality suppose that $j'_1\neq0$. Set
\begin{equation}\label{e1}
e_1=(1,0,\dots,0),\quad e_2=(0,1,\dots,0),\dots,\ e_d=(0,\dots,0,1),
\end{equation}
and $j''=j'-e_1$. Pick $\e<\min\{|\g'-\g|:\,\g\in\G\}$, and set
$$
\p_{\g',j'',\e}(y)=\frac{(y-\g')^{j''}}{j''!}\p\left(\frac{|y-\g'|}{\e}\right),
$$
where $\p$ is defined in (\ref{p0}). We have
 \begin{equation}\label{j}
\hat f(e^{2\pi i\la y,t\ra}\p_{\g',j'',\e}(y))=
 \sum_{\g\in\G}\sum_{\|j\|\le J}(-1)^{\|j\|}q_j(\g)D^j(e^{2\pi i\la y,t\ra}\p_{\g',j'',\e}(y))(\g)
\end{equation}
 Since
 $$
D^j(\p_{\g',j'',\e}(y))(\g)=\left\{\begin{array}{l} 0\mbox{ if }\g\neq\g'\mbox{ or }j\neq j'',
\\1\mbox{ if }\g=\g'\mbox{ and }j=j'',
  \end{array}\right.
$$
 we see that expression (\ref{j}) is equal to
 $$
 (-1)^Jq_{j'}(\g')2\pi it_1 e^{2\pi i\la\g',t\ra}+(-1)^J\sum_{s=2}^d q_{j''+e_s}(\g')2\pi it_s e^{2\pi i\la\g',t\ra}+(-1)^{J-1}q_{j''}(\g')e^{2\pi i\la\g',t\ra}.
 $$
 The first summand is unbounded in $t_1\in\R$, hence the function
  $$
  f(\hat\p_{\g',j'',\e}(x-t))=\hat f(e^{2\pi i\la y,t\ra}\p_{\g',j'',\e}(-y))
  $$
  is unbounded and not almost periodic. We obtain the contradiction, therefore, $J=0$ and $\hat f$ is a measure. \bs
\medskip

In order to prove Theorems \ref{T5} and \ref{T6}, we need  the following proposition:
\begin{Pro}\label{P4}
Let $\l_n,\tau_n\in\R^d$ be two sequences such that $\tau_n\to0$, $|\l_n|>|\l_{n-1}|+1$ for all $n$, and
\begin{equation}\label{g}
 \log|\tau_n|/\log|\l_n|\to-\infty\quad\mbox{as}\quad n\to\infty.
\end{equation}
Let $\mu$ be any measure from $S^*(\R^d)$ such that its restriction for each ball $B(\l_n,1/(2|\l_n|))$ equals $|\tau_n|^{-2/3}(\d_{\l_n+\tau_n}-\d_{\l_n})$.
Then there is $\psi\in S(\R^d)$  such that $\hat\mu\star\hat\psi(t)$ is unbounded. In particular, $\hat\mu$ is not  s-almost periodic  distribution.
\end{Pro}
{\bf Proof}.
By thinning out the sequence $\tau_n$, we can assume that for all $n$
\begin{equation}\label{s1}
\sum_{p<n}|\tau_p|^{-1/3}<(1/3)|\tau_n|^{-1/3},
\end{equation}
and
\begin{equation}\label{s2}
\sum_{p>n}|\tau_p|^{2/3}<(2/(3\pi))|\tau_n|^{2/3}.
\end{equation}
Set
 $$
 \psi(x)=\sum_n|\tau_n|^{1/3}\p(|\l_n||x-\l_n|),
 $$
 where $\p$ is defined in \eqref{p0}.  By \eqref{g}, $|\tau_n|=o(1/|\l_n|^T)$ as $n\to\infty$ for every $T<\infty$. Therefore, for all
  $k\in(\N\cup\{0\})^d,\ N\in\N$
we have $D^k\psi(x)=o(|\l_n|^{-N})$ for $x\in B(\l_n,1/(2|\l_n|))$. Hence, $(D^k\psi)(x)(1+|x|^N)$ is bounded on $\R^d$ for all $N$ and $k$, i.e., $\psi\in S(\R^d)$.
  By \eqref{p0}, $\psi(x)=0$ for $x\not\in\cup_n B(\l_n,1/(2\l_n))$. Hence, for every $t\in\R^d$
 $$
 \hat\mu(\hat\psi(t-y))=\mu(\psi(x)e^{-2\pi i\la x,t\ra})=\sum_{n=1}^\infty|\tau_n|^{-1/3}[\p(|\tau_n||\l_n|)e^{-2\pi i\la(\l_n+\tau_n),t\ra}-
 \p(0)e^{-2\pi i\la \l_n,t\ra}].
 $$
 For large $n$ we have $|\tau_n|<1/(3|\l_n|)$,  therefore, $\p(|\tau_n||\l_n|)=\p(0)=1$. Besides, for $t=\tau_n/(2|\tau_n|^2)$
 $$
 |e^{-2\pi i\la(\l_n+\tau_n),t\ra}-e^{-2\pi i\la \l_n,t\ra}|=|e^{-2\pi i\la\tau_n,t\ra}-1|=2.
  $$
  Therefore,
 \begin{equation}\label{s3}
  |\hat\mu(\hat\psi(t-y))|\ge2|\tau_n|^{-1/3}-2\sum_{p<n}|\tau_p|^{-1/3}-\sum_{p>n}|\tau_p|^{-1/3}|e^{-2\pi i\la\tau_p,t\ra}-1|.
 \end{equation}
Taking into account \eqref{s1}, \eqref{s2}, and the estmates
$$
|e^{-2\pi i\la\tau_p,t\ra}-1|\le 2\pi|\tau_p||t|=\pi|\tau_p||\tau_n|^{-1},
$$
we obtain that \eqref{s3} is more than $2|\tau_n|^{-1/3}/3$. Hence the convolution $(\hat\mu\star\hat\psi)(t)$ is unbounded
on the sequence $\tau_n/(2|\tau_n|^2)$, and the distribution  $\hat\mu$ is not s-almost periodic. \bs
\medskip

{\bf Proof of  Theorem \ref{T5}}. Suppose that $A$ is not $p$-discrete set. Then there are two sequences
$\l_n,\l'_n\in A$ such that $\l_n$ and $\tau_n:=\l'_n-\l_n$ satisfy \eqref{g} and $|\l_n|>1+|\l_{n-1}|$. Check that the  measure
$$
\mu=\sum_n|\tau_n|^{-2/3}[\d_{\l'_n}-\d_{\l_n}]
$$
belongs to $S^*(\R^d)$.

For any $\phi\in S(\R^d)$ we have
$$
  |(\mu,\phi)|  \le\sum_n |\tau_n|^{-2/3}|\phi(\l'_n)-\phi(\l_n)|  \le\sum_n |\tau_n|^{1/3}N_{0,1}(\phi),
$$
where $N_{0,1}(\phi)$ is defined in \eqref{d}. By \eqref{g}, $\tau_n=O(n^{-T})$ for any $T<\infty$, therefore the sum converges,  $\mu$ satisfies \eqref{d},
and $\mu$ is a temperate distribution. Applying Proposition \ref{P4}, we obtain that $\hat\mu$ is not s-almost periodic.

 Proposition \ref{P4}  actually implies  the unboundedness of the convolution $\hat\mu\star\hat\psi$ with some $\psi\in S(\R^d)$, which proves the last part of the theorem.

Sufficiency follows from Corollary \ref{C3}.  \bs
\medskip

Our proof of Theorem \ref{T6} uses the following lemmas:
\begin{Lem}[Y.Meyer \cite{M2}, Lemma 7]\label{L2}
Let $\a\in(0,1/6)$. For every integer $M>M_\a$ there exists an
$M$-periodic locally finite measure $\s=\s_M$ that is a sum of Dirac masses on
$\L_M =M^{-1}\Z\setminus[-\a M,\a M]$, and whose Fourier transform is also supported
by $\L_M$. To be precise,
$$
\supp\s=M^{-1}\Z\setminus\left[\cup_{k\in\Z}(kM+(-\a M,\a M)\right].
$$
\end{Lem}
For $\tau\in\R$ set $\s_M^\tau=\sum_{\l\in\supp\s_M}\d_{\l+\tau}$.
\begin{Lem}\label{L3}
Let $M\ge16,\,\a=1/8$. Then for any $\phi\in S(\R)$ and $\tau\in(0,1/2)$ we have
\begin{equation}\label{e}
 |(\s_M^{2\tau}-\s_M^\tau,\phi)|\le C'\tau N_{2,1}(\phi),
\end{equation}
where $N_{2,1}(\phi)$ is defined in \eqref{d}, and $C'$ is an absolute constant.
\end{Lem}
{\bf Proof of Lemma \ref{L3}}. Taking into account that $|\l|\ge\a M>2$ for every $\l\in\supp\s_M$, we get
$$
  |((\d_{\l+2\tau}-\d_{\l+\tau}),\phi)|\le\max_{0\le\t\le1}|\phi'(\l+(1+\t)\tau)|\tau\le\tau N_{2,1}(\phi)|\l|^{-2}.
  $$
Note that a number of points of $\supp\s_M^{2\tau}\cup\supp\s_M^\tau$ on the interval
$$
((k-1)M+(1/8)M,kM-(1/8)M),\,k\ge1,
$$
 is less than $2M^2$. Therefore,
 $$
|(\s_M^{2\tau}-\s_M^\tau,\phi)|\le 2M^2\tau N_{2,1}(\phi)\sum_{k=1}^\infty[(k-1)M+M/8]^{-2}.
$$
\bs

{\bf Proof of Theorem \ref{T6}}.  Let $\tau_n$ be a sequence of mutually independent
over $\Z$ positive numbers such that $\tau_n<1/16$ and $\log\tau_n/n\to-\infty$ as $n\to\infty$. Let $\L_{M_n}$ be the set defined in Lemma \ref{L2} with
$M_n=16^n,\ \a=1/8$, and $\s_{M_n}$ be the corresponding measure.
It is not hard to check that the sets $\L_{M_n}+\tau_n, \L_{M_p}+2\tau_p,\ n,p\in\N,$ are mutually disjoint if $n,\,p\ge n_0$ and $n_0$ large enough. Since
$\L_{M_p}\cap(-\a M_n,\a M_n)=\emptyset$ for $p\ge n$, we see that the measure
$$
  \nu=\sum_{n\ge n_0}\tau_n^{-2/3}(\s_{M_n}^{2\tau_n}-\s_{M_n}^{\tau_n})
$$
is locally finite for appropriate $n_0$. Using Lemma \ref{L3}, we obtain
$$
  |(\nu,\phi)|\le \sum_{n\ge n_0} C'\tau_n^{1/3}N_{2,1}(\phi).
$$
Since $\tau_n\le n^{-6}$ for $n$ large enough, we see that the sum converges, the measure $\nu$ satisfies condition \eqref{d}, and $\nu\in S^*(\R)$.
Moreover, its spectrum is a subset of the locally finite set $\cup_{n\ge n_0}\L_{M_n}$, therefore the measure $\mu=\check\nu$ is a crystalline measure.
By Proposition \ref{P4}, to check that $\mu$ is not s-almost periodic distribution we only need to find a sequence $\l_n\to\infty$ such that $\log\l_n=O(n)$ and
$$
(\l_n-1/(2\l_n),\l_n+1/(2\l_n))\cap\supp\nu=\{\l_n,\l_n+\tau_n\}.
$$
Fix $n\ge n_0$. Set $\eta_{n,j}=M_n+j/M_n$ with $j\in\N$ such that $M_n+M_n/8\le\eta_{n,j}\le 2M_n-M_n/8$,
and set $I_{n,j}:=(\eta_{n,j}-1/(2\eta_{n,j}),\eta_{n,j}+1/(2\eta_{n,j}))$.
 For every fixed $n$ these intervals do not intersect.
Since $2M_n-M_n/8+2\tau_n<M_p/8$ for $p>n$, we get
$$
I_{n,j}\cap[\supp\s_{M_p}^{2\tau_p}\cup\supp\s_{M_p}^{\tau_p}]=\emptyset \quad\forall j,\,\forall p>n.
$$
Then for every $p<n$ a number of points of the form $kM_p+q/M_p+\tau_p$ or $kM_p+q/M_p+2\tau_p$, $k,q\in\N\cup\{0\}$ on the interval $(M_n,2M_n)$
is at most $2M_nM_p$. Summing over all $p<n$, we get
$$
  \#\{\cup_{p<n}[\supp\s_{M_p}^{2\tau_p}\cup\supp\s_{M_p}^{\tau_p}]\}\cap(M_n,2M_n)<2M_n^2/15.
$$
On the other hand, the number of points $\eta_{n,j}$ on the interval $(M_n+M_n/8,\,2M_n-M_n/8)$ is $3M_n^2/4$, and the same is the number of intervals
$I_{n,j}\subset(M_n,2M_n)$.  Hence there is $j'$  such that
 $$
 I_{n,j'}\cap[\supp\s_{M_p}^{2\tau_p}\cup\supp\s_{M_p}^{\tau_p}]=\emptyset\quad\forall p\neq n.
$$
 Set $\l_n=\eta_{n,j'}+\tau_n$. Clearly, $(16)^n<\l_n<2(16)^n$. Since $\tau_n<4^{-1}16^{-n}$ for $n$ large enough, we see that
 $$
 \l_n,\, \l_n+\tau_n\in(\l_n-1/(2\l_n),\,\l_n+1/(2\l_n))\subset I_{n,j'},
 $$
 and
 $$
 \eta_{j,n}+\tau_n,\,\eta_{j,n}+2\tau_n\not\in(\l_n-1/(2\l_n),\,\l_n+1/(2\l_n))\quad\mbox{for}\quad j\neq j'.
 $$
 We obtain that the measure $\nu$ satisfies Proposition \ref{P4} for $n_0$ large enough.
 It follows from Corollary \ref{C2} that the measure $|\hat\nu|$ does not belong to $S^*(\R^d)$, hence $\hat\nu$ is not the Fourier quasicrystal. \bs

\section{Fourier coefficients and Fourier transform of almost periodic distributions}\label{S5}

For any almost periodic function $g$ on $\R^d$ its coefficient Fourier for an exponent $\l\in\R^d$ is defined by the formula
$$
   a(\l,g)=\lim_{R\to\infty}\frac{1}{\omega_d R^d}\int_{B(0,R)} g(t)e^{-2\pi i\la\l,t\ra}dt,
$$
where $\omega_d$ is the volume of the unit ball in $\R^d$, and the limit exists uniformly with respect to  shifts of $g$.

Now let $f$ be an almost periodic temperate distribution on $\R^d$. By L.Ronkin \cite{R},  its Fourier coefficient corresponding to an exponent $\l\in\R^d$  is defined
by the equality
$$
    a(\l,f)=a(\l,f\star\phi)/\hat\phi(\l),
$$
where $\phi$ is any $C^\infty$-function with compact support such that $\hat\phi(\l)\neq0$, and $a(\l,f\star\phi)$ is the Fourier coefficient
of the almost periodic function $f\star\phi$.
Note that  the definition matches with the previous one for an almost periodic function $g$. Indeed, in this case
$$
   a(\l,g\star\phi)=\lim_{R\to\infty}\frac{1}{\omega_d R^d}\int_{B(0,R)}\int_{\supp\phi} g(t-x)\phi(x) e^{-2\pi i\la\l,t\ra}dx\,dt
$$
$$
   = \int_{\supp\phi}\left[\lim_{R\to\infty}\frac{1}{\omega_d R^d}\int_{B(0,R)} g(t-x) e^{-2\pi i\la\l,(t-x)\ra}dt\right]\phi(x)e^{-2\pi i\la\l,x\ra}dx
   =a(\l,g)\hat\phi(\l).
$$
Then for any almost periodic temperate distribution $f$ and any $\phi,\,\psi\in C^\infty(\R^d)$ with compact supports such that $\hat\phi(\l)\neq0,\,\hat\psi(\l)\neq0$
we apply the above equality for almost periodic functions
$f\star\phi,\,f\star\psi$ and get
$$
   \frac{a(\l,f\star\phi)}{\hat\phi(\l)}=\frac{a(\l,(f\star\phi)\star\psi)}{\hat\phi(\l)\hat\psi(\l)}=\frac{a(\l,(f\star\psi)\star\phi))}{\hat\phi(\l)\hat\psi(\l)}
   =\frac{a(\l,f\star\psi)}{\hat\psi(\l)}.
$$
Therefore the definition of $a(\l,f)$ does not depend on a function $\phi$.

 Y.Meyer (\cite{M4}, Theorem 3.8) proved that for any almost periodic function $g$ the following conditions are equivalent

i) the sum $\sum_{|\l|<R}|a(\l,g)|$ converges for every $R<\infty$,

ii) $\hat g$ is a measure,

iii) $\hat g$ is an atomic measure of the form $\hat g=\sum_\l a(\l,g)\d_\l$.

Then Y.Meyer generalized this result to almost periodic measures.
\smallskip

The result has a simple generalization to almost periodic distributions and their Fourier coefficients in sense of Ronkin:
\begin{Th}\label{T7}
Let $f$ be an almost periodic distribution. Then the following conditions are equivalent:

i) the sum $\sum_{|\l|<R}|a(\l,f)|$ converges for every $R<\infty$,

ii) $\hat f$ is a measure,

iii) $\hat f$ is an atomic measure of the form $\hat f=\sum_\l a(\l,f)\d_\l$.
\end{Th}
{\bf Proof}. Check the implication i)$\Rightarrow$iii). Then for any function $\phi\in C^\infty$ with compact support the sum
$$
\sum_{|\l|<R} |a(\l,f\star\phi)|=\sum_{|\l|<R} |a(\l,f)||\hat\phi(\l)|
$$
converges for any $R<\infty$. Therefore $\widehat{f\star\phi}=\hat f\hat\phi$ is a measure of the form
$$
\sum_\l a(\l,f\star\phi)\d_\l=\sum_\l a(\l,f)\hat\phi(\l)\d_\l.
$$
Note that any function $\psi\in S(\R^d)$ can be approximated by $C^\infty$-functions with compact supports. Therefore for every $\psi\in S(\R^d)$
the sum
$$
\hat f\hat \psi=\sum_\l a(\l,f)\hat\psi(\l)\d_\l
$$
 is a continuous linear functional on the space of all continuous functions on every ball $\overline{B(0,R)},\,R<\infty$.
If we take $\psi$ such that $\hat\psi(\l)\equiv 1$ for $|\l|\le R$, we get that the restriction of $\hat f$ on $\overline{B(0,R)}$ coincides
with $\sum_\l a(\l,f)\d_\l$, and iii) is proved.

There is no need to prove the implication iii)$\Rightarrow$ii).

 In order to check the implication ii)$\Rightarrow$i), suppose that $\hat f$ is a measure and $\phi$ is any $C^\infty$-function with compact support. Then
  $\widehat{f\star\phi}=\hat f\hat\phi$ is also a measure. Since $f\star\phi$ is
 an almost periodic function, we see, by Meyer's Theorem, that the sum
$$
\sum_{|\l|<R} |a(\l,f)||\hat\phi(\l)|=\sum_{|\l|<R} |a(\l,f\star\phi)|
$$
 converges for any $R<\infty$. If $\hat\phi(\l)\neq0$ for all $\l\in\overline{B(0,R)}$, we obtain i).  \bs
\bigskip

 I want to thank the referee for a careful reading of my article and important remarks that forced the author to completely change Theorem 3.
I also would like to thank N.Lev, G.Reti for the inaccuracy they discovered in one important definition. Finally, I thank the Department of Mathematics and Computer Science
of the Jagiellonian University for their hospitality and Professor Lukasz Kosinski for his interest in my work and useful discussions.

\end{document}